\documentclass[10pt]{article}
\usepackage{amstext}
\usepackage{amsmath}
\usepackage{amsfonts}
\usepackage{amssymb}
\usepackage{amsthm}
\usepackage{hyperref}
\usepackage{url}
\usepackage{indentfirst}
\usepackage{tikz}
\usepackage{booktabs}
\usetikzlibrary{shapes.geometric}
\usepackage{filecontents}
\usepackage[noend]{algpseudocode}
\usepackage[margin=1.5in]{geometry}    
\usepackage{algorithm}
\usepackage{graphicx}
\usetikzlibrary{shapes,positioning}
\usepackage[utf8]{inputenc}
\usepackage{hyperref}
\usepackage{caption}

\usepackage{enumerate} 
\usepackage[justification=centering]{caption} 
\usepackage[export]{adjustbox}
\usepackage{caption, subcaption}

\newcommand{\cA}{\mathcal{A}}
\newcommand{\cB}{\mathcal{B}}

\newtheorem{theorem}{Theorem}[section]

\newtheorem{remark}{Remark}[section]

\usepackage{chngcntr}
\counterwithout{remark}{section}
\counterwithout{figure}{section}
\counterwithout{table}{section}
\counterwithout{theorem}{section}
\counterwithout{lemma}{section}
\counterwithout{proposition}{section}

\title{Lexicographically maximal edges of dual hypergraphs and Nash-solvability of tight game forms}

\author{Vladimir Gurvich\\
vgurvich@hse.ru and vladimir.gurvich@gmail.com\\
RUTCOR, Rutgers University, Piscataway, NJ, United States;\\
National Research University Higher School of Economics, Moscow, Russia
\and
Mariya Naumova\\
mnaumova@business.rutgers.edu\\
Rutgers Business School, Rutgers University, Piscataway, NJ, United States\\
}

\begin{document}
\maketitle
\begin{abstract}
Let   $\cA = \{A_1, \ldots, A_m\}$  and  $\cB = \{B_1, \ldots, B_n\}$  be 
a pair of dual multi-hypergraphs on the common ground set
$O = \{o_1, \ldots, o_k\}$. 
Note that each of  them may have embedded or equal edges.
An edge is called containment minimal 
(or just minimal, for short) if it is not a strict superset of another edge. 
Yet, equal minimal edges may exist. By duality,  
\begin{enumerate}[(i)]
    \item
    $A \cap B \neq \emptyset$  for  every pair  $A \in \cA$  and  $B \in \cB$;
    \item
 if  $A$  is minimal then for every  $o \in A$  there exists a $B \in \cB$  such that
$A \cap B = \{o\}$.

We will extend claim (ii) as follows. A linear order $\succ$  over $O$
 defines a unique lexicographic order $\succ_L$  over the   $2^O$.  
 Let  $A$  be a lexicographically maximal (lexmax) edge of  $\cA$. Then,
 \item
 $A$  is minimal and for every  $o \in A$  there exists a minimal $B \in \cB$  such that  
 $A \cap B = \{o\}$  and  $o \succeq o'$  for each  $o' \in B$.  
\end{enumerate}

This property has important applications in game theory implying 
Nash-solvability of tight game forms as shown in the old (1975 and 1989) work of the first author. 
Here we give a new, very short, proof of  (iii).  

Edges  $A$  and  $B$  mentioned in (iii)  can be found out in polynomial time. This is trivial if  $\cA$  and $\cB$  are given explicitly. 
Yet, it is true even if only $\cA$  is given, and not explicitly, but by a polynomial containment oracle, which for a subset $O_A \subseteq O$ answers in polynomial time  whether  $O_A$  contains an edge of  $\cA$.

{\bf AMS subjects: 91A05, 94D10, 06E30.}

{\bf Keywords:} Power set, dual multi-hypergraphs, lexicographical order, lexmax edge, Nash equilibrium, Nash-solvability, tight game form.
\end{abstract}

\section{Main concepts and results}  
\label{s0}
{\bf Lexicographical order over the subsets.}
A linear order  $\succ$  over  a set  $O = \{o_1, \ldots, o_k\}$  
uniquely determines a lexicographical order $\succ_L$
over the power set  $2^O$ of all $2^k$ subsets of  $O$  as follows.
Roughly speaking, the more small elements are avoided by a set - the larger it is.
In particular, $O' \succ_L  O''$  whenever  $O' \subset O''$ and, hence,
the empty set  $\emptyset \subset O$  is the largest in $2^O$.  

More precisely, to compare  two arbitrary subsets  $O', O'' \subseteq O$  
consider their symmetric difference
$\Delta = (O' \setminus O'') \cup (O'' \setminus O')$.
Clearly,  $\Delta \neq \emptyset$  if and only if  $O'$ and $O''$ are distinct.
Let  $o$  be the minimum with respect to  $\succ$  element in  $\Delta$.
If  $o \in   (O' \setminus O'')$ then  $O'' \succ_L O'$;
if  $o \in   (O'' \setminus O')$ then  $O' \succ_L O''$.

This definition can be equivalently reformulated as follows. 
Assume, without any loss of generality, that  $o_1 \prec \ldots  \prec o_k$,  
assign the negative weight  $w(o_i) = -2^{k-i}$  to every $o_i \in O$, 
and set   $w(U) = \sum_{o \in U} w(o)$   for each  $U \subseteq 0$. 
Then, $O' \succ_L O''$  if and only if  $w(O') > w(O'')$.

\bigskip

{\bf  Dual multi-hypergraphs.}  
We assume that the reader is familiar with basic notions
related to monotone Boolean functions, in particular, with DNFs and duality. 
An introduction can be found in \cite{CH11}; see Sections 1, 3 and 4.

\smallskip 

Multi-hypergraphs  $\cA = \{A_1, \ldots, A_m\}$  and  $\cB = \{B_1, \ldots, B_n\}$  
defined on the common ground set  $O = \{o_1, \ldots, o_k\}$  are called dual 
if  (i)  holds, 
$A \cap B \neq \emptyset$  for  every pair  $A \in \cA$  and  $B \in \cB$, 
and  also 

\smallskip 
(iv a) 
for each  $B^T \subseteq O$  such that 
$B^T \cap B \neq \emptyset$  for  every  $B \in \cB$ 
there exists an  $A \in \cA$  such that  $A \subseteq B^T$. 

\smallskip 

If (i) and (iv a)  hold we say that  $\cA$ is dual to $\cB$  and 
use notation $\cA = \cB^d$. 
Swapping  $\cA$  and  $\cB$  in (iv a) we obtain  (iv b) 
and an equivalent definition of duality, 
that is, (i)  and  (iv a)  hold  if and only  (i) and  (iv b)  hold.
In other words, $\cA = \cB^d$  if and only if  $\cB = \cA^d$. 
So we just say that  $\cA$  and $\cB$  are dual.

\medskip 

A multi-hypergraph is called {\em Sperner} 
if no two of its distinct edges contain one another;
in particular, they cannot be equal. 
In this case,
we have a hypergraph rather than multi-hypergraph. 
For a multi-hypergraph there exists a unique 
dual Sperner hypergraph. 
If  $\cA$  and  $\cB$  are dual and Sperner then 
$\cA^{dd} = \cA$  and  $\cB^{dd} = \cB$; 
furthermore, 
$\cup_{i=1}^m A_i = \cup_{j=1}^n B_j = O$. 
In general, for multi-hypergraphs, 
$\cup_{i=1}^m A_i$  and  $\cup_{j=1}^n B_j$  are 
subsets of  $O$  that may differ. 

\begin{remark} 
\label{r-Sperner} 
It is well known \cite{CH11} that 
(dual) multi-hypergraphs are in one-to-one correspondence with
(dual) monotone DNFs; 
(prime) implicants of the latter correspond to 
(minimal) edges of the former. 
Furthermore, Sperner hypergraphs correspond to irredundant DNFs.
However, we do not restrict ourselves to this case. 
Although our main result 
(the lexicographical theorem, see below) would not lose much but 
its applications to Nash-solvability would; 
see the last section.
\end{remark} 

Claims (i) and (ii) of the Abstract are well-known; see for example  \cite{CH11}. 
Actually, (i) is required directly by the definition of duality of $\cA$ and $\cB$. 
If (ii) fails then edge  $A$ cannot be minimal, since  
$A \setminus \{o\}$   still intersects all  $B \in \cB$.

Our main result is (iii). Fix an arbitrary order $\succ$ over $O$  and 
find a lexmax edge  $A^0$, that is, one  maximal  with respect to the lexicographical order  
$\succ_L$ over $2^A$.  
Note that  $A^0$  may be not unique, but all lexmax edges are equal.

\begin{theorem}
\label{t1}
A lexmax edge  $A^0 \in \cA$  is minimal.
Furthermore, for every  $o^0 \in A^0$  there exists a (minimal) edge  
$B^0 \in \cB$  such that  $A^0 \cap B^0 = \{o^0\}$  and  
$o^0 \succeq o$  for each   $o \in B^0$.
\end{theorem}

\proof  
A lexmax edge must be minimal,
since a set is strictly less than any its proper subset in order $\succ_L$.

Assume for contradiction that there exists an  $o^0 \in A^0$  such that
for every  minimal $B^0 \in \cB$  satisfying   $B^0 \cap A^0 = \{o^0\}$
there exists  an  $o \in B^0$  such that  $o \succ o^0$.  
Clearly, this assumption holds for every $B^0 \in \cB$ 
if it holds for each minimal  $B^0 \in \cB$.
Let us show that it contradicts the lexmaximality of  $A^0$.
To do so, partition all edges of  $\cB$ into two types:  
\begin{enumerate}[(a)]
    \item there exists an  $o \in  B \cap A^0$  distinct from  $o^0$.
    \item $B \cap A^0 = \{o^0\}$. 
\end{enumerate}
In case (b), by our assumption,  
there is an  $o \in B$  such that  $o \succ o^0$. 
In both cases, (a) and (b), choose the specified element $o$ 
from each  $B \in \cB$, thus, getting a transversal  $B^T$.
By (iv a),  there exists an  $A \in \cA$  such that
$A \subseteq B^T$ and, hence, $A \succeq_L B^T$.

Furthermore,  $B^T \succ_L A^0$.
Indeed, by construction,  $o^0 \not \in B^T$  and it is replaced  
by some larger elements, $o \succ o^0$, in case  (b), 
while all other elements of  $B^T$, if any, 
belong to  $A^0 \setminus \{o^0\}$, in accordance with case (a).

Thus, by transitivity,  $A \succ_L A^0$, 
while  $A^0$  is a lexmax edge of  $\cA$, by assumption of the theorem, 
which is a contradiction.
\qed

\section{Determining edges  $A^0$ and $B^0$  in polynomial time.}
\label{s2}
The problem is trivial when  multi-hypergraphs
$\cA$  and  $\cB$  are given explicitly.
We will solve it when only  $\cA$  is given, and not explicitly, 
but by a polynomial containment oracle. 
For an arbitrary subset  $O_A \subseteq O$  
this oracle answers in polynomial time the question 
$Q(\cA, O_A)$: whether  $O_A$  contains an edge $A \in \cA$.

By duality, 
$A \not \subseteq O_A$  for all  $A \in \cA$  if and only if
$O_B = O \setminus O_A$  contains an edge   $B \in \cB$.
In other words, $Q(\cA, O_A)$  is answered in the negative if and only if
$Q(\cB, O_B)$  is answered in the positive. Thus, we do not need
two separate oracles for  $\cA$ and $\cB$; 
it is sufficient to have one, say, for $\cA$.

\smallskip

{\bf Determining a lexmax edge $\cA^0$.}
Recall that multi-hypergraph  $\cA$  may contain several lexmax edges,
but they are all equal.

Fix an arbitrary linear order  $\succ$  over  $O$.
Without loss of generality we can assume that  
$o_1 \prec \ldots \prec o_k$.

\smallskip

{\bf Determining a lexmax edge $A^0$.} 
Recall that multi-hypergraph  $\cA$  may contain several lexmax edges  $A^0$,
but they are all equal. 
Fix an arbitrary linear order  $\succ$  over  $O$.
Wlog we can assume that  $O = \{o_1, \ldots, o_k\}$  
and  $o_1 \prec \ldots \prec o_k$.

\smallskip 

Step 1: 
Consider  $O^1_t = \{o_t, \ldots, o_k\}$ and, by asking question $Q(\cA, O^1_t)$ for  
$t = 1, \ldots, k$, find the maximum  $t_1$  for which the answer is still positive.
Then,  $o_{t_1}$  belongs to  $\cA^0$,  while  $o_1, \ldots, o_{t_1 - 1}$  do not.  

\smallskip 

Step 2: 
Consider  $O^2_t = \{o_{t_1}, o_{t_1 + t}, \ldots, o_k\}$  
and, by asking question $Q(\cA, O^2_t)$ for  
$t = 1, \ldots, k - t_1$, find the maximum  $t_2$  
for which the answer is still positive. 
Then,  $o_{t_1}, o_{t_1 + t_2} \in \cA^0$,  
while  $o_t \not \in \cA^0$  for any other  $t < t_1 + t_2$. 

\smallskip 

Step 3: 
Consider  $O^3_t = \{o_{t_1}, o_{t_1 + t_2}, o_{t_1 + t_2 + t}, \ldots, o_k\}$  
and, by asking question $Q(\cA, O^3_t)$ for  
$t = 1, \ldots, k - (t_1 + t_2)$, find the maximum  $t_3$  
for which the answer is still positive. 
Then,  $o_{t_1}, o_{t_1 + t_2}, o_{t_1 + t_2 + t_3} \in \cA^L$, 
while  $o_t \not \in \cA^0$  for any other  $t < t_1 + t_2 + t_3$;  etc. 

\smallskip 

This procedure will produce a lexmax edge $A^0$ 
in at most  $k$  polynomial iterations. 
Notice that on each step  $i$ 
we can speed up the search of $t_i$  by applying the dichotomy. 

\medskip 

{\bf Determining an edge  $B^0$  from Theorem \ref{t1}.}
First, find a lexmax edge  $A^0 \in \cA$  and choose an arbitrary  $o^0 \in A^0$. 
We look for an edge  $B^0 \in \cB$  such that  
$A^0 \cap B^0 = \{o^0\}$  and  $o^0 \succ o$  
for every  $o \in B^0 \setminus \{o^0\}$,  that is, 

$$B^0 \subseteq  O_B = 
O  \setminus [(A^0  \setminus \{o^0\}) \cup \{o \mid o \succ o^0\}].$$ 

By Theorem \ref{t1}, such  $B^0$  exists and, hence, 
the oracle answers  $Q(\cB, O_B)$ in the positive, 
or equivalently,  $Q(\cA, O \setminus O_B)$  in the negative. 
We could take any  $B^0 \in \cB$  such that  $B^0 \subseteq O_B$. 
Yet, multi-hypergraph  $\cB$  is not given explicitly. 
To get  $B^0$  we ``minimize"   $O_B$. 
To do so, we delete its elements one by one in some order until obtain  
a minimum set  $O^M$   for which the answer to 
 $Q(\cA, O \setminus O^M)$  is still negative, that is, 
to   $Q(\cA, O \setminus (O^M \setminus \{o\})$  
it becomes positive for all $o \in O^M$.  
Then, we set  $B^0 = O^M$. 

Note that again we can speed up the procedure by using dichotomy. Note also that this reduction procedure is not necessarily unique, since 
we can eliminate elements of  $O_B$ in any order.
Thus, in contrast to  $A^0$, several distinct edges  $B^0$  
may satisfy Theorem \ref{t1}.

\section{Applications to Nash equilibria}
\label{s3}
Alice and Bob  play  $m \times n$  normal form game defined as follows.
Their strategies are edges of two dual multi-hypergraphs
$\cA = \{A_1, \ldots, A_m\}$  and  $\cB = \{B_1, \ldots, B_n\}$. 
A pair of strategies  $(A \in \cA, B \in \cB)$  is called a {\em situation}.
To each situation $(A,B)$  assign an arbitrary outcome  $o \in A \cap B$.
By duality of $\cA$ and $\cB$,  the intersection is not empty.
The obtained mapping  $g : \cA \times \cB \rightarrow O$  
is called a {\em tight game form}.

For example, dual multi-hypergraphs 

$\cA = \{(o_1, o_2), (o_2, o_3),(o_2, o_3, o_5),(o_2, o_3),(o_3, o_4)\}$  and 

$\cB = \{(o_1, o_3), (o_2, o_3), (o_2, o_3, o_5), (o_2, o_3, o_6), (o_2,o_4), (o_2, o_4, o_6, o_7\}$

\medskip 
\noindent 
generate tight game form   

\medskip 
\begin{center}
\[\begin{bmatrix}
o_1 & o_2 & o_2 & o_2 & o_2 & o_2\\
o_3 & o_3 & o_2 & o_3 & o_2 & o_2\\
o_3 & o_2 & o_5 & o_2 & o_2 & o_2\\
o_3 & o_3 & o_2 & o_3 & o_2 & o_2\\
o_3 & o_3 & o_3 & o_3 & o_4 & o_4\\
\end{bmatrix}.\]
\end{center}
\medskip

A general (finite two-person) game form is defined as an arbitrary mapping
$g : X \times Y \rightarrow O$, where  $X, Y$, and  $O$  
are three arbitrary finite sets.
Duality of  $\cA$  and  $\cB$  outlines an important class of tight game forms.

\begin{remark}
\label{r2} 
We consider multi-hypergraphs rather than (Sperner) hypergraphs
to get a larger class of tight game forms and, thus, more applications.
\newline
Otherwise, it would be reasonable to restrict ourselves to Sperner hypergraphs. 
Theorem \ref{t1}  would not lose much, but its applications in game theory would.
Indeed, tight game forms generated by the dual Sperner hypergraphs 
form a very special subfamily in the family of all tight game forms, 
generated by dual multi-hypergraphs.
\end{remark}

Linear orders  $\succ_a$  and  $\succ_b$  over  $A$  and $B$
will be interpreted as preferences of Alice and Bob, respectively.
Both are maximizers. 
The triplet  $(g, \succ_a, \succ_b)$  is called
a {\em  two-person normal form game} or just a game, for short.

A situation  $(A,B)$  of game  $(g, \succ_a, \succ_b)$ 
is called a {\em Nash equilibrium} (NE) if
 
\smallskip

$g(A,B) \succeq_a g(A',B)$  for each  $A' \in \cA$  and 

\smallskip
  
$g(A,B) \succeq_b g(A,B')$  for each  $B' \in \cB$.

\smallskip
\noindent
In other words, a player cannot profit  
by changing her/his strategy provided the opponent
keeps his/her strategy unchanged; or differently, 
if $A$  is a best response to $B$ and  $B$  is a best response to  $A$.

We can interpret Theorem \ref{t1} as follows.
Let  $A^0$ be a lexmax strategy of Alice and 
$o^0$ be the best outcome of Bob in $A^0$. 
By Theorem \ref{t1}, there is a strategy  $B^0$  of Bob 
such that $A^0 \cap  B^0 = \{o^0\}$  and $o^0 \succeq_a o$  
for every $o \in B^0$.            
In other words, $(A^0, B^0)$  is a NE.
Indeed,  $B^0$  is a best response to  $A^0$, 
by definition of  $o^0$  and 
$A^0$  is a best response to  $B^0$, by Theorem \ref{t1}. 

This type of NE has an interesting property.
The strategy of Alice depends only on her preference $\succ_a$;
she may be unaware of Bob's preference, which is important for applications.
Such NE are called  {\em lexmax NE of Alice}.
{\em Lexmax NE of Bob} are defined similarly; 
see more details in  \cite{GN21}.

In case of games with opposite preferences 
these two types of  NE coincide.
Note that in this case the concept 
of lexmax is a refinement of maxmin and minmax. 
Lexmax NE of Alice and Bob may coincide in other cases too; 
for example, a NE may be just unique.

\medskip 

However, in general, if both players choose their lexmax strategies 
then the obtained situation may be not a NE.
Consider the set of outcomes 
$O = \{o_1, o_2, o_3\}$, the  $2 \times 2$  
game form  $g : X \times Y \rightarrow  O$  given by the table

\medskip 

$ (x_1, y_1) \; (x_1, y_2)  \;\;\;\;\;  o_1 \; o_2$

$ (x_2, y_1) \; (x_2, y_2)  \;\;\;\;\;  o_1 \; o_3$

\medskip
\noindent 
and preferences of Alice and Bob given by orders 
$\succ_A$  and  $\succ_B$   such that 
$o_2 \succ_A o_3$   and  
$o_2 \succ_B o_1 \succ_B o_3$. 
Note that place of $o_1$  in $\succ_A$  is irrelevant. 

In the obtained game 
$(g, \succ_A, \succ_B)$  the lexmax strategies of 
Alice and Bob are  $x_1$  and  $y_1$,  respectively.
Yet, $(x_1, y_1)$  is not a NE. 
Indeed, Bob can improve replacing $y_1$ by $y_2$ 
and getting $o_2$  instead of  $o_1$.
Furthermore, situations  $(x_1, y_2)$  and  $(x_2, y_1)$ 
are both NE, in accordance with Theorem \ref{t1}.

\bigskip

The existence of NE in two-person games with tight game forms
for arbitrary preferences of the players was shown
in \cite{Gur75}; 
see also \cite{Gur89} and \cite{GN21} for more details.
The case of opposite preferences 
was considered earlier by Edmonds and Fulkerson \cite{EF70};
see also \cite{Gur73}.

\section*{Acknowledgement}
This work was included in the research plan of
the Research University Higher School of Economics.
The authors are thankful to Endre Boros for helpful remarks and suggestions.


\begin{thebibliography}{99}
\bibitem{CH11}
Y. Crama and P. L. Hammer. Boolean functions: 
Theory, algorithms, and applications. Cambridge University Press, 2011.

\bibitem{EF70}
J. Edmonds and D.R. Fulkerson, Bottleneck extrema,
J. of Combinatorial Theory, 8 (1970) 299--306.

\bibitem{Gur73}
V. Gurvich, On theory of multi-step games, 
USSR Comput. Math. and Math. Phys. 13:6 (1973) 143--161.

\bibitem{Gur75}
V. Gurvich,
Solution of positional games in pure strategies, 
USSR Comput. Math. and Math. Phys. 15: 2 (1975) 74--87.

\bibitem{Gur89}
V. Gurvich,
Equilibrium in pure strategies,
Soviet Math. Dokl. 38:3 (1989) 597--602.

\bibitem{GN21}
V. Gurvich and M. Naumova,
Polynomial algorithms computing two lexicographically safe Nash equilibria 
in finite two-person games with tight game forms given by oracles. 
In: ArXiv e-prints (2021) https://arxiv.org/abs/2108.05469  

\end{thebibliography}
\end{document}